\documentclass[12pt]{article}%
\usepackage{amsfonts}
\usepackage{sw20bams}
\usepackage{amsmath}
\usepackage{amssymb}
\usepackage{graphicx}%
\providecommand{\U}[1]{\protect\rule{.1in}{.1in}}
\begin{document}

\title{Cantor-solus and Cantor-multus Distributions}
\author{Steven Finch}
\date{March 20, 2020}
\maketitle

\begin{abstract}
The Cantor distribution is obtained from bitstrings; the Cantor-solus
distribution (a new name) admits only strings without adjacent $1$ bits. \ We
review moments and order statistics associated with these. \ The Cantor-multus
distribution is introduced -- which instead admits only strings without
isolated $1$ bits -- and more complicated formulas emerge.

\end{abstract}

\footnotetext{Copyright \copyright \ 2020 by Steven R. Finch. All rights
reserved.}A bitstring is \textbf{solus} if all of its $1$s are isolated.
\ Such strings were called Fibonacci words (more fully, words obeying the
Fibonacci restriction) in \cite{Pr1-tcs3}. \ We shall reserve the name
Fibonacci for a different purpose, as in \cite{Fi1-tcs3, Fi2-tcs3}. \ 

A bitstring is \textbf{multus} if each of its $1$s possess at least one
neighboring $1$. \ Such strings were called good sequences in \cite{AG-tcs3}.
\ Counts of solus $n$-bitstrings have a quadratic character, whereas counts of
multus $n$-bitstrings have a cubic character. \ More on the meaning of this
and on other related combinatorics will appear later.

\section{Cantor Distribution}

Let $0<\vartheta\leq1/2$; for instance, we could take $\vartheta=1/3$ as in
the classical case. \ Let $\bar{\vartheta}=1-\vartheta$.\ \ Consider a mapping
\cite{LT-tcs3}
\[
F(\omega_{1}\omega_{2}\omega_{3}\cdots\omega_{m})=\frac{\bar{\vartheta}%
}{\vartheta}%
{\displaystyle\sum\limits_{i=1}^{m}}
\omega_{i}\vartheta^{i}%
\]
from the set $\Omega$ of finite bitstrings ($m<\infty$) to the nonnegative
reals. \ The $2^{m}$ bitstrings in $\Omega$ of length $m$ are assumed to be
equiprobable. \ Consider the generating function \cite{Pr1-tcs3}%
\[
G_{n}(z)=%
{\displaystyle\sum\limits_{\omega\in\Omega}}
F(\omega)^{n}z^{\left\vert \omega\right\vert }%
\]
where $\left\vert \omega\right\vert $ denotes the length of the bitstring.
\ Clearly%
\[
G_{0}(z)=%
{\displaystyle\sum\limits_{\omega\in\Omega}}
z^{\left\vert \omega\right\vert }=%
{\displaystyle\sum\limits_{m=0}^{\infty}}
2^{m}z^{m}=\frac{1}{1-2z}.
\]
The quantity%
\[
\frac{\lbrack z^{m}]G_{n}(z)}{[z^{m}]G_{0}(z)}=\frac{1}{2^{m}}[z^{m}]G_{n}(z)
\]
is the $n^{\text{th}}$ Cantor moment for strings of length $m$; let $\mu_{n}$
denote the limit of this as $m\rightarrow\infty$. \ Denote the empty string by
$\varepsilon$. \ From values%
\[%
\begin{array}
[c]{ccccc}%
F(\varepsilon)=0, &  & F(0\omega)=\vartheta\,F(\omega), &  & F(1\omega
)=\bar{\vartheta}+\vartheta\,F(\omega)
\end{array}
\]
and employing the recurrence \cite{SF-tcs3}%
\[
\Omega=\varepsilon+\{0,1\}\times\Omega,
\]
we have%
\begin{align*}
G_{n}(z)  &  =%
{\displaystyle\sum\limits_{\omega\in\Omega}}
\vartheta^{n}\,F(\omega)^{n}z^{1+\left\vert \omega\right\vert }+%
{\displaystyle\sum\limits_{\omega\in\Omega}}
\left(  \bar{\vartheta}+\vartheta\,F(\omega)\right)  ^{n}z^{1+\left\vert
\omega\right\vert }\\
&  =\vartheta^{n}zG_{n}(z)+z%
{\displaystyle\sum\limits_{i=0}^{n}}
\dbinom{n}{i}\bar{\vartheta}^{n-i}\vartheta^{i}G_{i}(z)\\
&  =2\vartheta^{n}zG_{n}(z)+z%
{\displaystyle\sum\limits_{i=0}^{n-1}}
\dbinom{n}{i}\bar{\vartheta}^{n-i}\vartheta^{i}G_{i}(z)
\end{align*}
for $n\geq1$; thus%
\[
G_{n}(z)=\frac{z}{1-2\vartheta^{n}z}%
{\displaystyle\sum\limits_{i=0}^{n-1}}
\dbinom{n}{i}\bar{\vartheta}^{n-i}\vartheta^{i}G_{i}(z).
\]
Dividing both sides by $G_{0}(z)$, we have \cite{LT-tcs3, Ev-tcs3, DF-tcs3,
DMRV-tcs3}%
\[
\mu_{n}=\frac{1}{2\left(  1-\vartheta^{n}\right)  }%
{\displaystyle\sum\limits_{i=0}^{n-1}}
\dbinom{n}{i}\bar{\vartheta}^{n-i}\vartheta^{i}\mu_{i}%
\]
because%
\[
\lim_{z\rightarrow z_{0}}\frac{z}{1-2\vartheta^{n}z}=\frac{1}{2\left(
1-\vartheta^{n}\right)  }%
\]
and the singularity $z_{0}=1/2$ of $G_{n}(z)$ is a simple pole. In particular,
when $\vartheta=1/3$,
\[%
\begin{array}
[c]{ccccc}%
\mu_{1}=1/2, &  & \mu_{2}=3/8, &  & \mu_{2}-\mu_{1}^{2}=1/8
\end{array}
\]
and, up to small periodic fluctuations \cite{DMRV-tcs3, GW-tcs3, GP-tcs3},%
\[
\mu_{n}\sim C\,n^{-\ln(2)/\ln(3)},
\]%
\[
C=\frac{1}{2\ln(3)}%
{\displaystyle\int\limits_{0}^{\infty}}
\left(
{\displaystyle\prod\limits_{k=2}^{\infty}}
\frac{1+e^{-2x/3^{k}}}{2}\right)  e^{-2x/3}x^{\ln(2)/\ln(3)-1}dx=0.733874...
\]
as $n\rightarrow\infty$.

We merely mention a problem involving order statistics. \ Let $\xi_{n}$ denote
the expected value of the minimum of $n$ independent Cantor-distributed random
variables. \ It is known that \cite{Ho-tcs3}
\[
\xi_{n}=\frac{1}{2^{n}-2\vartheta}\left[  \bar{\vartheta}+\vartheta%
{\displaystyle\sum\limits_{i=1}^{n-1}}
\dbinom{n}{i}\xi_{i}\right]
\]
in general. \ In the special case $\vartheta=1/3$, it follows that%
\[%
\begin{array}
[c]{ccccccccc}%
\xi_{1}=1/2, &  & \xi_{2}=3/10, &  & \xi_{3}=1/5, &  & \xi_{4}=33/230, &  &
\xi_{5}=5/46
\end{array}
\]
and, up to small periodic fluctuations \cite{KP-tcs3},%
\[
\xi_{n}\sim c\,n^{-\ln(3)/\ln(2)},
\]%
\[
c=\frac{2}{3\ln(2)}\Gamma\left(  \frac{\ln(3)}{\ln(2)}\right)  \zeta\left(
\frac{\ln(3)}{\ln(2)}\right)  =1.9967049717...
\]
as $n\rightarrow\infty$. \ If $\eta_{n}$ denotes the expected value of the
maximum of $n$ variables, then%
\[
1-\eta_{n}\sim c\,n^{-\ln(3)/\ln(2)}%
\]
by symmetry.

A final problem concerns the sum of all moments of the classical Cantor
distribution \cite{Pr2-tcs3}:%
\begin{align*}%
{\displaystyle\sum\limits_{n=0}^{\infty}}
\mu_{n}  &  =-\frac{1}{3}+\frac{2}{3}%
{\displaystyle\sum\limits_{k=1}^{\infty}}
\left(  \frac{2}{3}\right)  ^{k}%
{\displaystyle\sum\limits_{j=1}^{2^{k}}}
\frac{1}{j}\\
&  =3.3646507281...
\end{align*}
answering a question asked in \cite{DRAK-tcs3}.

\section{Cantor-solus Distribution}

We examine here the set $\Omega$ of finite solus bitstrings ($m<\infty$).
\ Let%
\[%
\begin{array}
[c]{ccccc}%
f_{k}=f_{k-1}+f_{k-2}, &  & f_{0}=0, &  & f_{1}=1
\end{array}
\]
denote the Fibonacci numbers. \ The $f_{m+2}$ bitstrings in $\Omega$ of length
$m$ are assumed to be equiprobable. \ Clearly%
\[
G_{0}(z)=%
{\displaystyle\sum\limits_{\omega\in\Omega}}
z^{\left\vert \omega\right\vert }=%
{\displaystyle\sum\limits_{m=0}^{\infty}}
f_{m+2}z^{m}=\frac{1+z}{1-z-z^{2}}.
\]
From additional values%
\[%
\begin{array}
[c]{ccc}%
F(1)=\bar{\vartheta}, &  & F(10\omega)=\bar{\vartheta}+\vartheta^{2}F(\omega)
\end{array}
\]
and employing the recurrence \cite{SF-tcs3}%
\[
\Omega=\varepsilon+1+\{0,10\}\times\Omega,
\]
we have%
\begin{align*}
G_{n}(z)  &  =\bar{\vartheta}^{n}z+%
{\displaystyle\sum\limits_{\omega\in\Omega}}
\vartheta^{n}\,F(\omega)^{n}z^{1+\left\vert \omega\right\vert }+%
{\displaystyle\sum\limits_{\omega\in\Omega}}
\left(  \bar{\vartheta}+\vartheta^{2}F(\omega)\right)  ^{n}z^{2+\left\vert
\omega\right\vert }\\
&  =\bar{\vartheta}^{n}z+\vartheta^{n}zG_{n}(z)+z^{2}%
{\displaystyle\sum\limits_{i+j=n}}
\dbinom{n}{i,j}\bar{\vartheta}^{i}\vartheta^{2j}G_{j}(z)\\
&  =\bar{\vartheta}^{n}z+\vartheta^{n}zG_{n}(z)+\vartheta^{2n}z^{2}%
G_{n}(z)+z^{2}%
{\displaystyle\sum\limits_{\substack{i+j=n,\\j<n}}}
\dbinom{n}{i,j}\bar{\vartheta}^{i}\vartheta^{2j}G_{j}(z)
\end{align*}
for $n\geq1$; thus%
\[
G_{n}(z)=\frac{1}{1-\vartheta^{n}z-\vartheta^{2n}z^{2}}\left[  \bar{\vartheta
}^{n}z+z^{2}%
{\displaystyle\sum\limits_{\substack{i+j=n,\\j<n}}}
\dbinom{n}{i,j}\bar{\vartheta}^{i}\vartheta^{2j}G_{j}(z)\right]  .
\]
The purpose of using multinomial coefficients here, rather than binomial
coefficients as in Section 1, is simply to establish precedent for Section 3.
\ Let $\varphi=(1+\sqrt{5})/2=1.6180339887...$ be the Golden mean. \ Dividing
both sides by $G_{0}(z)$, we have \cite{Pr1-tcs3}%
\begin{align*}
\mu_{n}  &  =\frac{1}{1-\vartheta^{n}/\varphi-\vartheta^{2n}/\varphi^{2}%
}\left[  0+\frac{1}{\varphi^{2}}%
{\displaystyle\sum\limits_{\substack{i+j=n,\\j<n}}}
\dbinom{n}{i,j}\bar{\vartheta}^{i}\vartheta^{2j}\mu_{j}\right] \\
&  =\frac{1}{\varphi^{2}-\vartheta^{n}\varphi-\vartheta^{2n}}%
{\displaystyle\sum\limits_{\substack{i+j=n,\\j<n}}}
\dbinom{n}{i,j}\bar{\vartheta}^{i}\vartheta^{2j}\mu_{j}%
\end{align*}
because%
\[
\lim_{z\rightarrow z_{0}}\frac{\bar{\vartheta}^{n}z}{G_{0}(z)}=\lim
_{z\rightarrow z_{0}}\frac{1-z-z^{2}}{1+z}\bar{\vartheta}^{n}z=0
\]
and the singularity $z_{0}=1/\varphi$ of $G_{n}(z)$ is a simple pole. \ In
particular, when $\vartheta=1/3$,
\[%
\begin{array}
[c]{ccccc}%
\mu_{1}=0.338826..., &  & \mu_{2}=0.203899..., &  & \mu_{2}-\mu_{1}%
^{2}=0.089096...
\end{array}
\]
and, up to small periodic fluctuations,%
\[
\mu_{n}\sim(0.616005...)n^{-\ln(\varphi)/\ln(3)}(3/4)^{n},
\]
as $n\rightarrow\infty$. \ An integral formula in \cite{Pr1-tcs3} for the
preceding numerical coefficient involves a generating function of exponential
type:%
\[
M(x)=e^{-x/3}%
{\displaystyle\sum\limits_{k=0}^{\infty}}
\frac{\mu_{k}}{k!}\left(  \frac{4x}{9}\right)  ^{k},
\]
namely%
\[
\frac{1}{2\varphi\ln(3)}%
{\displaystyle\int\limits_{0}^{\infty}}
M(x)e^{-2x/3}x^{\ln(\varphi)/\ln(3)-1}dx
\]
(we believe that the fifth decimal given in \cite{Pr1-tcs3} is incorrect,
perhaps a typo). Unlike the formula for $C$ earlier, this expression depends
on the sequence $\mu_{1}$, $\mu_{2}$, $\mu_{3}$, \ldots\ explicitly.

With regard to order statistics, it is known that \cite{CP-tcs3}%
\[
\xi_{n}=\frac{1}{1-\vartheta\varphi^{-n}-\vartheta^{2}\varphi^{-2n}}\left[
\bar{\vartheta}\varphi^{-2n}+\vartheta%
{\displaystyle\sum\limits_{i=1}^{n-1}}
\dbinom{n}{i}\varphi^{-i}\varphi^{-2(n-i)}\xi_{i}\right]  ,
\]%
\[
\eta_{n}=\frac{1}{1-\vartheta\varphi^{-n}-\vartheta^{2}\varphi^{-2n}}\left[
\bar{\vartheta}\left(  1-\varphi^{-n}\right)  +\vartheta^{2}%
{\displaystyle\sum\limits_{j=1}^{n-1}}
\dbinom{n}{j}\varphi^{-2j}\varphi^{-(n-j)}\eta_{j}\right]
\]
in general. \ In the special case $\vartheta=1/3$, we have, up to small
periodic fluctuations,%
\[
\xi_{n}\sim(3.31661...)n^{-\ln(3)/\ln(\varphi)},
\]%
\[
3/4-\eta_{n}\sim(5.35114...)n^{-\ln(3)/\ln(\varphi)}%
\]
as $n\rightarrow\infty$.

\section{Cantor-multus Distribution}

We examine here the set $\Omega$ of finite multus bitstrings ($m<\infty$).
\ Let%
\[%
\begin{array}
[c]{ccccc}%
f_{k}=2f_{k-1}-f_{k-2}+f_{k-3}, &  & f_{0}=0, &  & f_{1}=f_{2}=1
\end{array}
\]
denote the second upper Fibonacci numbers \cite{Krc-tcs3}. \ The $f_{m+2}$
bitstrings in $\Omega$ of length $m$ are assumed to be equiprobable. \ Clearly%
\[
G_{0}(z)=%
{\displaystyle\sum\limits_{\omega\in\Omega}}
z^{\left\vert \omega\right\vert }=%
{\displaystyle\sum\limits_{m=0}^{\infty}}
f_{m+2}z^{m}=\frac{1-z+z^{2}}{1-2z+z^{2}-z^{3}}.
\]
From additional values%
\[
F(11\omega)=\bar{\vartheta}+\bar{\vartheta}\vartheta+\vartheta^{2}F(\omega),
\]%
\[
F(1110\omega)=\bar{\vartheta}+\bar{\vartheta}\vartheta+\bar{\vartheta
}\vartheta^{2}+\vartheta^{4}F(\omega)
\]
and employing the recurrence%
\[
\Omega=\varepsilon+1+\{0,11,1110\}\times\Omega,
\]
we have%
\begin{align*}
G_{n}(z)  &  =\bar{\vartheta}^{n}z+%
{\displaystyle\sum\limits_{\omega\in\Omega}}
\vartheta^{n}\,F(\omega)^{n}z^{1+\left\vert \omega\right\vert }+%
{\displaystyle\sum\limits_{\omega\in\Omega}}
\left(  \bar{\vartheta}+\bar{\vartheta}\vartheta+\vartheta^{2}F(\omega
)\right)  ^{n}z^{2+\left\vert \omega\right\vert }\\
&  +%
{\displaystyle\sum\limits_{\omega\in\Omega}}
\left(  \bar{\vartheta}+\bar{\vartheta}\vartheta+\bar{\vartheta}\vartheta
^{2}+\vartheta^{4}F(\omega)\right)  ^{n}z^{4+\left\vert \omega\right\vert }\\
&  =\bar{\vartheta}^{n}z+\vartheta^{n}zG_{n}(z)+z^{2}%
{\displaystyle\sum\limits_{i+j+k=n}}
\dbinom{n}{i,j,k}\bar{\vartheta}^{i}\left(  \bar{\vartheta}\vartheta\right)
^{j}\left(  \vartheta^{2}\right)  ^{k}G_{k}(z)\\
&  +z^{4}%
{\displaystyle\sum\limits_{i+j+k+\ell=n}}
\dbinom{n}{i,j,k,\ell}\bar{\vartheta}^{i}\left(  \bar{\vartheta}%
\vartheta\right)  ^{j}\left(  \bar{\vartheta}\vartheta^{2}\right)  ^{k}\left(
\vartheta^{4}\right)  ^{\ell}G_{\ell}(z)\\
&  =\bar{\vartheta}^{n}z+\vartheta^{n}zG_{n}(z)+\vartheta^{2n}z^{2}%
G_{n}(z)+z^{2}%
{\displaystyle\sum\limits_{\substack{i+j+k=n,\\k<n}}}
\dbinom{n}{i,j,k}\bar{\vartheta}^{i+j}\vartheta^{j+2k}G_{k}(z)\\
&  +\vartheta^{4n}z^{4}G_{n}(z)+z^{4}%
{\displaystyle\sum\limits_{\substack{i+j+k+\ell=n,\\\ell<n}}}
\dbinom{n}{i,j,k,\ell}\bar{\vartheta}^{i+j+k}\vartheta^{j+2k+4\ell}G_{\ell}(z)
\end{align*}
for $n\geq1$; thus%
\begin{align*}
G_{n}(z)  &  =\frac{1}{1-\vartheta^{n}z-\vartheta^{2n}z^{2}-\vartheta
^{4n}z^{4}}\left[  \bar{\vartheta}^{n}z+z^{2}%
{\displaystyle\sum\limits_{\substack{i+j+k=n,\\k<n}}}
\dbinom{n}{i,j,k}\bar{\vartheta}^{i+j}\vartheta^{j+2k}G_{k}(z)\right. \\
&  \left.  +z^{4}%
{\displaystyle\sum\limits_{\substack{i+j+k+\ell=n,\\\ell<n}}}
\dbinom{n}{i,j,k,\ell}\bar{\vartheta}^{i+j+k}\vartheta^{j+2k+4\ell}G_{\ell
}(z)\right]  .
\end{align*}
Let%
\[
\psi=\frac{1}{3}\left[  2+\left(  \frac{25+3\sqrt{69}}{2}\right)
^{1/3}+\left(  \frac{25-3\sqrt{69}}{2}\right)  ^{1/3}\right]
=1.7548776662...
\]
be the second upper Golden mean \cite{Krc-tcs3, Fi3-tcs3}. \ Dividing both
sides by $G_{0}(z)$, we have%
\begin{align*}
\mu_{n}  &  =\frac{1}{1-\vartheta^{n}/\psi-\vartheta^{2n}/\psi^{2}%
-\vartheta^{4n}/\psi^{4}}\left[  0+\frac{1}{\psi^{2}}%
{\displaystyle\sum\limits_{\substack{i+j+k=n,\\k<n}}}
\dbinom{n}{i,j,k}\bar{\vartheta}^{i+j}\vartheta^{j+2k}\mu_{k}\right. \\
&  \left.  +\frac{1}{\psi^{4}}%
{\displaystyle\sum\limits_{\substack{i+j+k+\ell=n,\\\ell<n}}}
\dbinom{n}{i,j,k,\ell}\bar{\vartheta}^{i+j+k}\vartheta^{j+2k+4\ell}\mu_{\ell
}\right] \\
&  =\frac{1}{\psi^{4}-\vartheta^{n}\psi^{3}-\vartheta^{2n}\psi^{2}%
-\vartheta^{4n}}\left[  \psi^{2}%
{\displaystyle\sum\limits_{\substack{i+j+k=n,\\k<n}}}
\dbinom{n}{i,j,k}\bar{\vartheta}^{i+j}\vartheta^{j+2k}\mu_{k}\right. \\
&  \left.  +%
{\displaystyle\sum\limits_{\substack{i+j+k+\ell=n,\\\ell<n}}}
\dbinom{n}{i,j,k,\ell}\bar{\vartheta}^{i+j+k}\vartheta^{j+2k+4\ell}\mu_{\ell
}\right]
\end{align*}
because%
\[
\lim_{z\rightarrow z_{0}}\frac{\bar{\vartheta}^{n}z}{G_{0}(z)}=\lim
_{z\rightarrow z_{0}}\frac{1-2z+z^{2}-z^{3}}{1-z+z^{2}}\bar{\vartheta}^{n}z=0
\]
and the singularity $z_{0}=1/\psi$ of $G_{n}(z)$ is a simple pole. \ In
particular, when $\vartheta=1/3$,
\[%
\begin{array}
[c]{ccccc}%
\mu_{1}=0.504968..., &  & \mu_{2}=0.416013..., &  & \mu_{2}-\mu_{1}%
^{2}=0.161020...
\end{array}
\]
but no asymptotics for $\mu_{n}$ are known. \ Order statistics likewise remain open.

\section{Bitsums}

We turn to a more fundamental topic: given a set $\Omega$ of finite
bitstrings, what can be said about the bitsum $S_{n}$ of a random $\omega
\in\Omega$ of length $n$? \ If $\Omega$ is unconstrained, i.e., if all $2^{n}$
strings are included in the sample, then%
\[%
\begin{array}
[c]{ccc}%
\mathbb{E}(S_{n})=n/2, &  & \mathbb{V}(S_{n})=n/4
\end{array}
\]
because a sum of $n$ independent Bernoulli($1/2$) variables is Binomial($n$%
,$1/2$). \ Expressed differently, the average density of $1$s in a random
unconstrained string is $1/2$, with a corresponding variance $1/4$. \ 

Let us impose constraints. \ If $\Omega$ consists of solus bitstrings, then
the total bitsum $a_{n}$ of all $\omega\in\Omega$ of length $n$ has generating
function \cite{Sn1-tcs3, GPF-tcs3}%
\[%
{\displaystyle\sum\limits_{n=0}^{\infty}}
a_{n}z^{n}=\frac{z}{\left(  1-z-z^{2}\right)  ^{2}}=z+2z^{2}+5z^{3}%
+10z^{4}+20z^{5}+\cdots
\]
and the total bitsum squared $b_{n}$ has generating function%
\[%
{\displaystyle\sum\limits_{n=0}^{\infty}}
b_{n}z^{n}=\frac{z\left(  1-z+z^{2}\right)  }{\left(  1-z-z^{2}\right)  ^{3}%
}=z+2z^{2}+7z^{3}+16z^{4}+38z^{5}+\cdots;
\]
hence $c_{n}=f_{n+2}b_{n}-a_{n}^{2}$ has generating function
\[%
{\displaystyle\sum\limits_{n=0}^{\infty}}
c_{n}z^{n}=\frac{z\left(  1-z\right)  }{(1+z)^{3}\left(  1-3z+z^{2}\right)
^{2}}=z+2z^{2}+10z^{3}+28z^{4}+94z^{5}+\cdots
\]
where $f_{n}$ is as in Section 2. \ Standard techniques \cite{SF-tcs3} give
asymptotics \
\[
\lim_{n\rightarrow\infty}\frac{\mathbb{E}(S_{n})}{n}=\lim_{n\rightarrow\infty
}\frac{a_{n}}{nf_{n+2}}=\frac{5-\sqrt{5}}{10}=0.2763932022...,
\]%
\[
\lim_{n\rightarrow\infty}\frac{\mathbb{V}(S_{n})}{n}=\lim_{n\rightarrow\infty
}\frac{c_{n}}{nf_{n+2}^{2}}=\frac{1}{5\sqrt{5}}=0.0894427190...
\]
for the average density of $1$s in a random solus string and corresponding
variance. \ 

If instead $\Omega$ consists of multus bitstrings, then the total bitsum
$a_{n}$ of all $\omega\in\Omega$ of length $n$ has generating function
\cite{Sn2-tcs3}%
\[%
{\displaystyle\sum\limits_{n=0}^{\infty}}
a_{n}z^{n}=\frac{z^{2}(2-z)}{\left(  1-2z+z^{2}-z^{3}\right)  ^{2}}%
=2z^{2}+7z^{3}+16z^{4}+34z^{5}+\cdots
\]
and the total bitsum squared $b_{n}$ has generating function%
\[%
{\displaystyle\sum\limits_{n=0}^{\infty}}
b_{n}z^{n}=\frac{z^{2}\left(  4-7z+4z^{2}+3z^{3}-z^{4}\right)  }{\left(
1-2z+z^{2}-z^{3}\right)  ^{3}}=4z^{2}+17z^{3}+46z^{4}+116z^{5}+\cdots;
\]
hence $c_{n}=f_{n+2}b_{n}-a_{n}^{2}$ has generating function
\[%
{\displaystyle\sum\limits_{n=0}^{\infty}}
c_{n}z^{n}=\frac{z^{2}\left(  4-9z+9z^{2}-9z^{3}-6z^{4}+z^{5}-6z^{6}%
+z^{8}\right)  }{\left(  1-z+2z^{2}-z^{3}\right)  ^{3}\left(  1-2z-3z^{2}%
-z^{3}\right)  ^{2}}=4z^{2}+19z^{3}+66z^{4}+236z^{5}+\cdots
\]
where $f_{n}$ is as in Section 3. \ We obtain asymptotics \
\begin{align*}
\lim_{n\rightarrow\infty}\frac{\mathbb{E}(S_{n})}{n}  &  =\lim_{n\rightarrow
\infty}\frac{a_{n}}{nf_{n+2}}\\
&  =\frac{1}{3}\left[  2-\left(  \frac{23+3\sqrt{69}}{1058}\right)
^{1/3}+\left(  \frac{-23+3\sqrt{69}}{1058}\right)  ^{1/3}\right] \\
&  =0.5885044113...,
\end{align*}%
\begin{align*}
\lim_{n\rightarrow\infty}\frac{\mathbb{V}(S_{n})}{n}  &  =\lim_{n\rightarrow
\infty}\frac{c_{n}}{nf_{n+2}^{2}}\\
&  =\frac{1}{1587}\left(  \frac{69}{2}\right)  ^{1/3}\left[  \left(
404685+35053\sqrt{69}\right)  ^{1/3}+\left(  404685-35053\sqrt{69}\right)
^{1/3}\right] \\
&  =0.2810976123...
\end{align*}
for the average density of $1$s in a random multus string and corresponding
variance. \ Unsurprisingly $0.588>1/2>0.276$ and $0.281>1/4>0.089$; a clumping
of $1$s forces a higher density than a separating of $1$s.

A famous example of an infinite aperiodic solus bitstring is the Fibonacci
word \cite{Fi1-tcs3, Fi2-tcs3}, which is the limit obtained recursively
starting with $0$ and satisfying substitution rules $0\mapsto01$, $1\mapsto0$.
\ The density of $1$s in this word is $1-1/\varphi\approx0.382$
\cite{Gry-tcs3}, which exceeds the average $0.276$ but falls well within the
one-sigma upper limit $0.276+\sqrt{0.089}=0.575$. \ We wonder if an
analogously simple construction might give an infinite aperiodic multus
bitstring with known density.

\section{Longest Bitruns}

We turn to a different topic:\ given a set $\Omega$ of finite bitstrings, what
can be said about the duration $R_{n,1}$ of the longest run of $1$s in a
random $\omega\in\Omega$ of length $n$? \ If $\Omega$ is unconstrained, then
\cite{SF-tcs3}%
\[
\mathbb{E}(R_{n,1})=\frac{1}{2^{n}}\left[  z^{n}\right]
{\displaystyle\sum\limits_{k=1}^{\infty}}
\left(  \frac{1}{1-2z}-\frac{1-z^{k}}{1-2z+z^{k+1}}\right)  ,
\]
the Taylor expansion of the numerator series is \cite{Sn3-tcs3}%
\[
z+4z^{2}+11z^{3}+27z^{4}+62z^{5}+138z^{6}+300z^{7}+643z^{8}+1363z^{9}%
+2866z^{10}+\cdots
\]
and, up to small periodic fluctuations \cite{Byd-tcs3, Sch-tcs3},%
\[
\mathbb{E}(R_{n,1})\sim\frac{\ln(n)}{\ln(2)}-\left(  \frac{3}{2}-\frac{\gamma
}{\ln(2)}\right)
\]
as $n\rightarrow\infty$. \ Of course, identical results hold for $R_{n,0}$,
the duration of the longest run of $0$s in $\omega$. \ 

If $\Omega$ consists of solus bitstrings, then it makes little sense to talk
about $1$-runs. \ For $0$-runs, over all $\omega\in\Omega$, we have%
\[
\mathbb{E}(R_{n,0})=\frac{1}{f_{n+2}}\left[  z^{n}\right]
{\displaystyle\sum\limits_{k=1}^{\infty}}
\left(  \frac{1+z}{1-z-z^{2}}-\frac{1+z-z^{k}-z^{k+1}}{1-z-z^{2}+z^{k+1}%
}\right)
\]
and the Taylor expansion of the numerator series is \cite{Sn3-tcs3}
\[
z+4z^{2}+9z^{3}+18z^{4}+34z^{5}+62z^{6}+110z^{7}+192z^{8}+331z^{9}%
+565z^{10}+\cdots
\]
where $f_{n}$ is as in Section 2.

If instead $\Omega$ consists of multus bitstrings, then we can talk both about
$1$-runs \cite{Sn3-tcs3}:%
\[
\mathbb{E}(R_{n,1})=\frac{1}{f_{n+2}}\left[  z^{n}\right]  \left\{  \frac
{-z}{\left(  1-z\right)  \left(  1-z+z^{2}\right)  }+%
{\displaystyle\sum\limits_{k=1}^{\infty}}
\left(  \frac{1+z^{2}}{1-2z+z^{2}-z^{3}}-\frac{1+z^{2}-z^{k-1}-z^{k}%
}{1-2z+z^{2}-z^{3}+z^{k+1}}\right)  z\right\}  ,
\]%
\[
\operatorname{num}=2z^{2}+7z^{3}+16z^{4}+32z^{5}+62z^{6}+118z^{7}%
+221z^{8}+409z^{9}+751z^{10}+\cdots
\]
and $0$-runs:%
\[
\mathbb{E}(R_{n,0})=\frac{1}{f_{n+2}}\left[  z^{n}\right]
{\displaystyle\sum\limits_{k=1}^{\infty}}
\left(  \frac{1+z^{2}}{1-2z+z^{2}-z^{3}}-\frac{1+z^{2}-z^{k-1}+z^{k}-2z^{k+1}%
}{1-2z+z^{2}-z^{3}+z^{k+2}}\right)  z,
\]%
\[
\operatorname{num}=z+2z^{2}+5z^{3}+11z^{4}+23z^{5}+45z^{6}+87z^{7}%
+165z^{8}+309z^{9}+573z^{10}+\cdots
\]
where $f_{n}$ is as in Section 3. \ Proof: the number of multus bitstrings
with no runs of $k$ $1$s has generating function \cite{Sn4-tcs3}%
\[%
\begin{array}
[c]{ccccccc}%
\dfrac{1+z^{2}-z^{k-1}-z^{k}}{1-2z+z^{2}-z^{3}+z^{k+1}}z &  & \text{if }k>1; &
& \dfrac{z}{1-z} &  & \text{if }k=1;
\end{array}
\]
we conclude by use of the summation identity
\[%
{\displaystyle\sum\limits_{j=0}^{\infty}}
j\cdot h_{j}(z)=%
{\displaystyle\sum\limits_{k=0}^{\infty}}
\left(
{\displaystyle\sum\limits_{i=0}^{\infty}}
h_{i}(z)-%
{\displaystyle\sum\limits_{i=0}^{k}}
h_{i}(z)\right)  .
\]
Study of runs of $k$ $0$s proceeds analogously \cite{Sn5-tcs3}. The solus and
multus results here are new, as far as is known. \ Asymptotics would be good
to see someday. \ 

\section{Acknowledgements}

I\ am thankful to Alois Heinz for helpful discussions and for providing the
generating function associated with $4,$ $19,$ $66,$ $236,$ $\ldots$ via the
Maple \textit{gfun} package; R\ and Mathematica have been useful throughout. I
am also indebted to a friend, who wishes to remain anonymous, for giving
encouragement and\ support (in these dark days of the novel coronavirus outbreak).

\end{document}